# ON LERCH'S TRANSCENDENT AND THE GAUSSIAN RANDOM WALK

By A. J. E. M. Janssen and J. S. H. van Leeuwaarden[1]

*Philips Research and EURANDOM*

Let $X_1, X_2, \ldots$ be independent variables, each having a normal distribution with negative mean $-\beta < 0$ and variance 1. We consider the partial sums $S_n = X_1 + \cdots + X_n$, with $S_0 = 0$, and refer to the process $\{S_n : n \geq 0\}$ as the Gaussian random walk. We present explicit expressions for the mean and variance of the maximum $M = \max\{S_n : n \geq 0\}$. These expressions are in terms of Taylor series about $\beta = 0$ with coefficients that involve the Riemann zeta function. Our results extend Kingman's first-order approximation [*Proc. Symp. on Congestion Theory* (1965) 137–169] of the mean for $\beta \downarrow 0$. We build upon the work of Chang and Peres [*Ann. Probab.* **25** (1997) 787–802], and use Bateman's formulas on Lerch's transcendent and Euler–Maclaurin summation as key ingredients.

**1. Introduction.** Let $X_1, X_2, \ldots$ be independent variables, each having a normal distribution with mean $-\beta < 0$ and variance 1. We consider the partial sums $S_n = X_1 + \cdots + X_n$, with $S_0 = 0$, and refer to the process $\{S_n : n \geq 0\}$ as the Gaussian random walk. In this paper we present explicit expressions for several characteristics of the distribution of the maximum

$$(1.1) \qquad M = \max\{S_n : n \geq 0\}.$$

The distribution of $M$ plays an important role in several areas of applied probability. In queueing theory, it typically occurs in a regime called *heavy traffic* (see [2, 18, 20, 26]), in which the load is just below its critical level, and so the queue is only just stable with relatively large queue lengths and waiting times. For the limiting waiting time $W = \lim_{n \to \infty} W_n$, with $W_1 = 0$ and $W_{n+1} = (W_n + X_n)^+$ (with $x^+ := \max\{0, x\}$), it follows from Spitzer's random-walk identities that $W$ is, in distribution, equal to $M$. In

Received January 2006; revised September 2006.

[1]Supported by the Netherlands Organisation for Scientific Research (NWO).

*AMS 2000 subject classifications.* 11M06, 30B40, 60G50, 60G51, 65B15.

*Key words and phrases.* Gaussian random walk, all-time maximum, Lerch's transcendent, Riemann zeta function, Spitzer's identity, Euler–Maclaurin summation.

the context of queues and heavy traffic, Kingman [20] was the first to observe the relevance of $M$ in his 1965 paper. He noticed among other things:

> Despite the apparent simplicity of the problem, there does not seem to be an explicit expression even for $\mathbb{E}M$..., but it is possible to give quite sharp inequalities and asymptotic results for small $\beta$.

Indeed, Kingman showed that, for $\beta \downarrow 0$,

$$(1.2) \qquad \mathbb{E}M = \frac{1}{2\beta} - c + \mathcal{O}(\beta), \qquad c \approx 0.58.$$

Determining the tail distribution of $M$ is tantamount to computing level crossing probabilities of the Gaussian random walk, that is, for $x > 0$, $\{M > x\} = \{\tau(x) < \infty\}$, where $\tau(x) = \inf\{n \geq 1 : S_n > x\}$. This level-crossings interpretation makes that the tail distribution of $M$ is important in *sequential analysis* and *risk theory*. Chang and Peres [10] derived an exact expression (2.1) for the expected value of the first descending ladder height (actually, they consider the first ascending ladder height for the Gaussian random walk with positive drift), denoted as $\mathbb{E}S_{\tau_-}$, with $\tau_- = \inf\{n \geq 1 : S_n \leq 0\}$, which by the relation $\mathbb{E}S_{\tau_-} = -\beta/\mathbb{P}(M = 0)$ (see [2], page 225) leads to an exact expression for $\mathbb{P}(M = 0)$. They present $\mathbb{E}S_{\tau_-}$ as a Taylor series about $\beta = 0$ with coefficients that involve the Riemann zeta function, a considerable achievement that generalizes first-order approximations of Siegmund [22] and second-order approximations of Chang [9].

Ladder heights fulfill an important role in probability theory, both in the exact analysis of random walks (see [2, 14]), and in the asymptotic analysis of boundary crossing problems [23]. In the latter case, a quantity of interest is the limiting expected overshoot, defined as $\mathbb{E}(S_\tau^2)/(2\mathbb{E}S_\tau)$, $\tau = \tau(0)$, for $\beta = 0$. This quantity can be shown to be $-\zeta(1/2)/\sqrt{2\pi} \approx 0.5826$, with $\zeta(z)$ the Riemann zeta function. The same quantity arises in sequentially testing for the drift of a Brownian motion [11], corrected diffusion approximations [22], simulation of Brownian motion [3, 8], option pricing [6] and thermodynamics of a polymer chain [12]. These applications have in common that a Brownian motion is observed only at equidistant sampling points. As it turns out (Kingman [20]) presents $c$ as $(2\pi)^{-1/2} \sum_{n=1}^{\infty} [\sqrt{n}(\sqrt{n} + \sqrt{n-1})^2]^{-1/2}$, which by Euler–Maclaurin summation can be shown to be $-(2\pi)^{-1/2}\zeta(1/2)$. Similar relations are the topic of Problem 602 posed by Glasser and Boersma in [15]), the $c$ in (1.2) is in fact $-\zeta(1/2)/\sqrt{2\pi}$, so Kingman, albeit in disguised form, related $\mathbb{E}M$ to the Riemann zeta function already in 1965. We shall extend Kingman's approximation (1.2) to an explicit expression for $\mathbb{E}M$, in the same spirit as Chang and Peres extended the results of Siegmund [22] and Chang [9]. Moreover, we present a similar expression for the variance of $M$, to be denoted by $\text{Var}\,M$. The new expressions for $\mathbb{E}M$ and $\text{Var}\,M$ both concern Taylor series about zero with coefficients that involve the Riemann zeta function.

The maximum and the first ladder height have been studied in the general setting of random walks with generally distributed increments (see [4, 9, 23]), the Gaussian random walk being a special case. For this general setting, Taylor series for the expected first ladder height and the expected maximum are presented in [4]. On a formal level, the results of Blanchet and Glynn [4] generalize our results and those of Chang and Peres [10]. However, finding the coefficients of the formal description of the Taylor series in [4] is a nontrivial exercise and requires the expansion of a characteristic function and the numerical evaluation of an integral (see [4], Section 6, in which an outline for this numerical procedure is given). This does not lead to exact expressions for the coefficients as in [10] or as in the present paper.

We first derive the Chang and Peres result [see (2.1) below] in our own fashion. Like Chang and Peres, we start from a Spitzer-type expression for $\mathbb{P}(M = 0)$, take its derivative with respect to $\beta$, rewrite the derivative in terms of the Riemann zeta function, and finally integrate to obtain (2.1). For rewriting the derivative, Chang and Peres built upon the 1905 paper of Hardy [17] and present an analytic continuation of the function $\mathrm{Li}_s(z) = \sum_{n=1}^{\infty} n^{-s} z^n$, known as the *polylogarithm* or *Jonquières function*. They were probably unaware of the fact that $\mathrm{Li}_s(z)$ is a special case of *Lerch's transcendent* [see (2.4)], for which the matter of analytic continuation has been established in full generality by Bateman (and/or the staff of the Bateman Manuscript Project); see [13], Section 1.11(8) and (2.5). Hence, although Chang and Peres [10] give a separate proof, their Theorem 2.1 should be attributed to Bateman.

Our derivation of (2.1)—that incorporates Bateman's formulas and an asymptotic determination of the integration constant—sets the stage for the derivation of the new explicit expressions for $\mathbb{E}M$ and $\mathrm{Var}\,M$. As an aside, we obtain the following asymptotic results for $\beta \downarrow 0$:

$$(1.3) \qquad \mathbb{E}M = \frac{1}{2\beta} + \frac{\zeta(1/2)}{\sqrt{2\pi}} + \frac{1}{4}\beta + \mathcal{O}(\beta^2)$$

and

$$(1.4) \qquad \mathrm{Var}\,M = \frac{1}{4\beta^2} - \frac{1}{4} - \frac{2\zeta(-1/2)}{\sqrt{2\pi}}\beta - \frac{1}{24}\beta^2 + \mathcal{O}(\beta^3),$$

where $\zeta(1/2) \approx -1.4604$ and $\zeta(-1/2) \approx -0.2079$. In comparing (1.2) and (1.3), (1.3) contains an additional term $\frac{1}{4}\beta$. This term, and $-\frac{1}{24}\beta^2$ in (1.4), follow from a rather intricate application of the Euler–Maclaurin summation formula. The error terms in both (1.3) and (1.4) will be replaced by Taylor series with coefficients that involve the Riemann zeta function.

1.1. *Structure of the paper.* We present our main results in the next section. Section 3 is devoted to an exposition of our derivation of the Chang and Peres result. The proofs of the new expressions for the mean and variance of the maximum are given in Sections 4 and 5, respectively. The new expressions for the mean and variance of $M$ are alternatives for their Spitzer-type counterparts. The latter tend to converge more slowly for a decreasing drift $\beta$, whereas the opposite holds for the new expressions. We investigate this difference in speed of convergence in Section 6. Concluding remarks are made in Section 7.

**2. Main results.** We present three theorems. The first, on $\mathbb{P}(M = 0)$, is essentially due to Chang and Peres [10], but we give a separate proof in Section 3:

THEOREM 1 (Chang and Peres [10]). *The probability that the maximum of the Gaussian random walk is zero satisfies*

$$(2.1) \qquad \mathbb{P}(M = 0) = \sqrt{2}\beta \exp\left\{ \frac{\beta}{\sqrt{2\pi}} \sum_{r=0}^{\infty} \frac{\zeta(1/2 - r)}{r!(2r+1)} \left( \frac{-\beta^2}{2} \right)^r \right\},$$

*for $0 < \beta < 2\sqrt{\pi}$.*

Then, largely motivated by Chang and Peres, but taking our own approach, we prove the next two theorems.

THEOREM 2. *The expectation of the maximum of the Gaussian random walk satisfies*

$$(2.2) \quad \mathbb{E}M = \frac{1}{2\beta} + \frac{\zeta(1/2)}{\sqrt{2\pi}} + \frac{1}{4}\beta + \frac{\beta^2}{\sqrt{2\pi}} \sum_{r=0}^{\infty} \frac{\zeta(-1/2 - r)}{r!(2r+1)(2r+2)} \left( \frac{-\beta^2}{2} \right)^r,$$

*for $0 < \beta < 2\sqrt{\pi}$.*

THEOREM 3. *The variance of the maximum of the Gaussian random walk satisfies*

$$(2.3) \quad \begin{aligned} \mathrm{Var}\, M &= \frac{1}{4\beta^2} - \frac{1}{4} - \frac{2\zeta(-1/2)}{\sqrt{2\pi}}\beta - \frac{\beta^2}{24} \\ &\quad - \frac{2\beta^3}{\sqrt{2\pi}} \sum_{r=0}^{\infty} \frac{\zeta(-3/2 - r)}{r!(2r+1)(2r+2)(2r+3)} \left( \frac{-\beta^2}{2} \right)^r, \end{aligned}$$

*for $0 < \beta < 2\sqrt{\pi}$.*

The key ingredients for obtaining the above series are Euler–Maclaurin summation and a result on Lerch's transcendent. Lerch's transcendent is defined as the analytic continuation of the series

$$(2.4) \qquad \Phi(z, s, v) = \sum_{n=0}^{\infty} (v + n)^{-s} z^n,$$

which converges for any real number $v \neq 0, -1, -2, \ldots$ if $z$ and $s$ are any complex numbers with either $|z| < 1$, or $|z| = 1$ and $\text{Re}(s) > 1$. Note that $\zeta(s) := \Phi(1, s, 1)$. We shall use the important result derived by Bateman [13], Section 1.11(8) [with $\zeta(s, v) := \Phi(1, s, v)$ the Hurwitz zeta function]:

$$(2.5) \qquad \Phi(z, s, v) = \frac{\Gamma(1-s)}{z^v} (\ln 1/z)^{s-1} + z^{-v} \sum_{r=0}^{\infty} \zeta(s-r, v) \frac{(\ln z)^r}{r!},$$

which holds for $|\ln z| < 2\pi$, $s \neq 1, 2, 3, \ldots$, and $v \neq 0, -1, -2, \ldots$.

**3. Proof of Theorem 1.** From Spitzer's identity for random walks [24] we have

$$(3.1) \quad \mathbb{P}(M = 0) = \exp\left\{ -\sum_{n=1}^{\infty} \frac{1}{n} \mathbb{P}(S_n > 0) \right\} = \exp\left\{ -\sum_{n=1}^{\infty} \frac{1}{n} P(-\beta\sqrt{n}) \right\},$$

with $P(\cdot)$ the standard normal distribution function

$$(3.2) \qquad P(a) = \frac{1}{\sqrt{2\pi}} \int_{-\infty}^{a} e^{-(1/2)x^2} \, dx.$$

The second equality in (3.1) follows from the normality of $S_n$.

With $F$ defined by

$$(3.3) \qquad F(\beta) = \sum_{n=1}^{\infty} \frac{1}{n} \frac{1}{\sqrt{2\pi}} \int_{\beta\sqrt{n}}^{\infty} e^{-(1/2)x^2} \, dx, \qquad \beta > 0,$$

we have

$$(3.4) \quad \begin{aligned} F'(\beta) &= \frac{-1}{\sqrt{2\pi}} \sum_{n=1}^{\infty} \frac{e^{-(1/2)\beta^2 n}}{\sqrt{n}} = \frac{-e^{-(1/2)\beta^2}}{\sqrt{2\pi}} \sum_{n=0}^{\infty} \frac{e^{-(1/2)\beta^2 n}}{\sqrt{n+1}} \\ &= \frac{-e^{-(1/2)\beta^2}}{\sqrt{2\pi}} \Phi\left(z = e^{-(1/2)\beta^2}, s = \frac{1}{2}, v = 1\right). \end{aligned}$$

Then by (2.5), when $0 < \frac{1}{2}\beta^2 < 2\pi$,

$$F'(\beta) = \frac{-e^{-(1/2)\beta^2}}{\sqrt{2\pi}} \left[ \frac{\Gamma(1/2)}{e^{-(1/2)\beta^2}} \left(\frac{1}{2}\beta^2\right)^{-1/2} \right.$$

$$(3.5) \qquad\qquad\qquad + e^{(1/2)\beta^2} \sum_{r=0}^{\infty} \zeta\left(\frac{1}{2} - r\right) \frac{(-(1/2)\beta^2)^r}{r!} \Bigg]$$

$$= \frac{-1}{\sqrt{2}} \left(\frac{1}{2}\beta^2\right)^{-1/2} - \frac{1}{\sqrt{2\pi}} \sum_{r=0}^{\infty} \zeta\left(\frac{1}{2} - r\right) \frac{(-(1/2)\beta^2)^r}{r!},$$

with $\zeta(s)$ denoting the Riemann zeta function. Thus, we get

$$(3.6) \quad F'(\beta) + \frac{1}{\beta} = \frac{-1}{\sqrt{2\pi}} \sum_{r=0}^{\infty} \zeta\left(\frac{1}{2} - r\right) \frac{(-(1/2)\beta^2)^r}{r!}, \qquad 0 < \beta < 2\sqrt{\pi}.$$

The series on the right-hand side of (3.6) converges uniformly in $\beta \in [0, \beta_0]$ when $0 \le \beta_0 < 2\sqrt{\pi}$; see (6.3). Therefore, when we integrate the identity in (3.6) from 0 to $\beta < 2\sqrt{\pi}$, we may interchange the sum and integral at the right-hand side, and we get

$$(3.7) \qquad F(\beta) + \ln\beta = L - \frac{1}{\sqrt{2\pi}} \sum_{r=0}^{\infty} \frac{\zeta(1/2 - r)(-1/2)^r \beta^{2r+1}}{r!(2r+1)},$$

where $L = \lim_{\beta \downarrow 0}(F(\beta) + \ln\beta)$.

We shall show that $L = -\frac{1}{2}\ln 2$. To that end, we note that

$$F(\beta) = \sum_{n=1}^{\infty} \frac{1}{n} \frac{1}{\sqrt{\pi}} \int_{\sqrt{(1/2)\beta^2 n}}^{\infty} e^{-u^2} \, du$$

$$= \frac{1}{2} \sum_{n=1}^{\infty} \frac{1}{n} \left(\frac{2}{\sqrt{\pi}} \int_{\sqrt{(1/2)\beta^2 n}}^{\infty} e^{-u^2} \, du - e^{-(1/2)\beta^2 n}\right)$$

$$(3.8) \qquad - \frac{1}{2} \ln(1 - e^{-(1/2)\beta^2})$$

$$= \frac{1}{4}\beta^2 \sum_{n=1}^{\infty} \frac{1}{(1/2)\beta^2 n} \left(\frac{2}{\sqrt{\pi}} \int_{\sqrt{(1/2)\beta^2 n}}^{\infty} e^{-u^2} \, du - e^{-(1/2)\beta^2 n}\right)$$

$$- \ln\beta + \frac{1}{2}\ln 2 + o(1),$$

as $\beta \downarrow 0$. The function

$$(3.9) \qquad g(y) := \frac{1}{y}\left(\frac{2}{\sqrt{\pi}} \int_{\sqrt{y}}^{\infty} e^{-u^2} \, du - e^{-y}\right), \qquad y > 0,$$

decays exponentially as $y \to \infty$, while $g(y) = \mathcal{O}(y^{-1/2})$, $y \downarrow 0$. It is then routine to show that

$$(3.10) \qquad \frac{1}{4}\beta^2 \sum_{n=0}^{\infty} g(\tfrac{1}{2}\beta^2 n) \to \frac{1}{2} \int_{0}^{\infty} g(y) \, dy, \qquad \tfrac{1}{2}\beta^2 \downarrow 0.$$

The latter integral can be evaluated as

$$
\int_0^\infty \frac{1}{y}\left( \frac{2}{\sqrt{\pi}} \int_{\sqrt{y}}^\infty e^{-u^2}\, du - e^{-y} \right) dy
$$

$$
= \left( \frac{2}{\sqrt{\pi}} \int_{\sqrt{y}}^\infty e^{-u^2}\, du - e^{-y} \right) \ln y \Big|_0^\infty
$$

(3.11)
$$
- \int_0^\infty \left( \frac{2}{\sqrt{\pi}} \cdot -\frac{1}{2} y^{-1/2} \cdot e^{-y} + e^{-y} \right) \ln y\, dy
$$

$$
= \frac{1}{\sqrt{\pi}} \int_0^\infty y^{-1/2} e^{-y} \ln y\, dy - \int_0^\infty e^{-y} \ln y\, dy
$$

$$
= \frac{1}{\sqrt{\pi}} \Gamma'(1/2) - \Gamma'(1) = -2\ln 2,
$$

by Abramowitz and Stegun ([1], 6.3.1-4, 258). Hence, $L = -\frac{1}{2}\ln 2$ indeed, and so it is shown that, for $0 < \beta < 2\sqrt{\pi}$, we have

$$
(3.12) \quad F(\beta) = -\ln \beta - \frac{1}{2}\ln 2 - \frac{1}{\sqrt{2\pi}} \sum_{r=0}^\infty \frac{\zeta(1/2 - r)(-1/2)^r \beta^{2r+1}}{r!(2r+1)},
$$

which, by (3.1), completes the proof of Theorem 1.

To recapitulate, we started from the Spitzer-type expression (3.1), rewrote its derivative (3.3) in terms of Lerch's transcendent (3.4), applied Bateman's formulas to obtain a Taylor series (3.5), integrated the Taylor series (3.7), and finally determined the integration constant $L$.

REMARK 4. The integration constant could have been determined from the relation $\mathbb{P}(M = 0) = -\beta/\mathbb{E}S_{\tau_-}$ (with $\mathbb{E}S_{\tau_-}$ the expected value of the first descending ladder height; see Section 1) and using the fact that $\mathbb{E}S_{\tau_-} = -1/\sqrt{2}$ for $\beta = 0$, as proven by Spitzer [25]; see also [21]. Alternatively, one could use the first-order approximation in [18], that is, $\mathbb{P}(M = 0) = \sqrt{2}\beta(1 + o(1))$ as $\beta \downarrow 0$. The primary purpose of this section, however, is to set the stage for the next two sections, in which there is no other way of determining integration constants than to apply asymptotic methods.

4. **Proof of Theorem 2.** From Spitzer's identity [24], we know that

$$
(4.1) \quad \mathbb{E}M = \sum_{n=1}^\infty \frac{1}{n} \mathbb{E}(S_n^+) = \sum_{n=1}^\infty \left( \frac{e^{-(1/2)\beta^2 n}}{\sqrt{2\pi n}} - \beta P(-\beta\sqrt{n}) \right).
$$

We then have

$$
\sum_{n=1}^\infty \frac{e^{-(1/2)\beta^2 n}}{\sqrt{2\pi n}} = \frac{e^{-(1/2)\beta^2}}{\sqrt{2\pi}} \Phi\left( z = e^{-(1/2)\beta^2}, s = \frac{1}{2}, v = 1 \right)
$$

(4.2)
$$= \frac{1}{\beta} + \frac{1}{\sqrt{2\pi}} \sum_{r=0}^{\infty} \frac{\zeta(1/2-r)}{r!} \left(-\frac{1}{2}\beta^2\right)^r.$$

Now we consider

(4.3)
$$G(\beta) = \sum_{n=1}^{\infty} \frac{1}{\sqrt{2\pi}} \int_{\beta\sqrt{n}}^{\infty} e^{-x^2/2} \, dx.$$

We have

(4.4)
$$G'(\beta) = \sum_{n=1}^{\infty} \frac{1}{\sqrt{2\pi}} \cdot \sqrt{n} \cdot - e^{-(1/2)\beta^2 n}$$

$$= \frac{-e^{-(1/2)\beta^2}}{\sqrt{2\pi}} \Phi\left(z = e^{-(1/2)\beta^2}, s = -\frac{1}{2}, v = 1\right).$$

Then by (2.5), when $\frac{1}{2}\beta^2 < 2\pi$,

(4.5)
$$G'(\beta) = \frac{-e^{-(1/2)\beta^2}}{\sqrt{2\pi}}$$

$$\times \left[ \frac{\Gamma(3/2)}{e^{-(1/2)\beta^2}} \left(\frac{1}{2}\beta^2\right)^{-3/2} + e^{(1/2)\beta^2} \sum_{r=0}^{\infty} \zeta\left(-\frac{1}{2}-r\right) \frac{(-(1/2)\beta^2)^r}{r!} \right]$$

$$= \frac{-1}{2\sqrt{2}} \left(\frac{1}{2}\beta^2\right)^{-3/2} - \frac{1}{\sqrt{2\pi}} \sum_{r=0}^{\infty} \zeta\left(-\frac{1}{2}-r\right) \frac{(-(1/2)\beta^2)^r}{r!}.$$

Therefore, we get

(4.6)
$$G'(\beta) = -\beta^{-3} - H(\beta);$$

$$H(\beta) = \frac{1}{\sqrt{2\pi}} \sum_{r=0}^{\infty} \frac{\zeta(-1/2-r)(-1/2)^r}{r!} \beta^{2r}.$$

We note that $H(\beta)$ is well behaved in $0 \le \beta < 2\sqrt{\pi}$, and that

(4.7)
$$\frac{d}{d\beta}\left[ G(\beta) - \frac{1}{2\beta^2} \right] = G'(\beta) + \frac{1}{\beta^3} = -H(\beta).$$

By integration from 0 to $\beta$, we thus get

(4.8)
$$G(\beta) - \frac{1}{2\beta^2} - \lim_{\varepsilon \downarrow 0} \left( G(\varepsilon) - \frac{1}{2\varepsilon^2} \right) = -\int_0^\beta H(\beta_1) \, d\beta_1$$

$$= \frac{-1}{\sqrt{2\pi}} \sum_{r=0}^{\infty} \frac{\zeta(-1/2-r)(-1/2)^r}{r!(2r+1)} \beta^{2r+1}.$$

We shall show that

$$(4.9) \qquad \lim_{\varepsilon \downarrow 0} \left( G(\varepsilon) - \frac{1}{2\varepsilon^2} \right) = -\frac{1}{4}.$$

To that end, we use the Euler–Maclaurin summation formula (see [7], Section 3.6, pages 40–42)

$$
\begin{aligned}
(4.10) \qquad \sum_{n=1}^{N} f(n) &= \int_1^N f(x)\,dx + \frac{1}{2}f(1) + \frac{1}{2}f(N) \\
&\quad + \sum_{k=1}^{m} \frac{B_{2k}}{(2k)!} \left( f^{(2k-1)}(N) - f^{(2k-1)}(1) \right) \\
&\quad - \int_1^N f^{(2m)}(x) \frac{B_{2m}(x - \lfloor x \rfloor)}{(2m)!}\,dx,
\end{aligned}
$$

where the $B_n(t)$ denote the Bernoulli polynomials, defined by

$$(4.11) \qquad \frac{z e^{zt}}{e^z - 1} = \sum_{n=0}^{\infty} \frac{B_n(t) z^n}{n!},$$

and the $B_n = B_n(0)$ denote the Bernoulli numbers. We apply (4.10) for $m = 1$, $N \to \infty$ and

$$(4.12) \qquad f_\delta(x) = \frac{1}{\sqrt{\pi}} \int_{\sqrt{\delta x}}^{\infty} e^{-u^2}\,du =: g(\delta x); \qquad \delta = \frac{1}{2}\varepsilon^2.$$

Hence,

$$
\begin{aligned}
(4.13) \qquad G_N(\varepsilon) &= \sum_{n=1}^{N} \frac{1}{\sqrt{2\pi}} \int_{\varepsilon\sqrt{n}}^{\infty} e^{-x^2/2}\,dx = \sum_{n=1}^{N} f_\delta(n) \\
&= \int_1^N g(\delta x)\,dx + \frac{1}{2}g(\delta) + \frac{1}{2}g(N\delta) \\
&\quad + \frac{1}{2}B_2(g'(N\delta) - g'(\delta))\delta - \int_1^N \delta^2 g''(\delta x) \frac{B_2(x - \lfloor x \rfloor)}{2}\,dx.
\end{aligned}
$$

Letting $N \to \infty$ and noting that for $g(y) = \frac{1}{\sqrt{\pi}} \int_{\sqrt{y}}^{\infty} e^{-u^2}\,du$ there holds that $g$, $g'$, $g'' \to 0$ exponentially fast as $y \to \infty$, we get

$$
\begin{aligned}
(4.14) \qquad G(\varepsilon) &= \int_1^\infty g(\delta x)\,dx + \frac{1}{2}g(\delta) - \frac{1}{2}B_2 g'(\delta)\delta \\
&\quad - \frac{1}{2}\int_1^\infty \delta^2 g''(\delta x) \frac{B_2(x - \lfloor x \rfloor)}{2}\,dx.
\end{aligned}
$$

Since $|B_{2k}(x)| \leq B_{2k}$ for $0 \leq x \leq 1$ (see [1], 23.1.13, page 805), the last integral at the right-hand side of (4.14) can be bounded by

$$(4.15) \qquad \int_1^\infty \delta^2 |g''(\delta x)| \tfrac{1}{2} B_2 \, dx = \tfrac{1}{12} \delta \int_\delta^\infty |g''(y)| \, dy.$$

We further get

$$(4.16) \qquad g'(y) = -\frac{e^{-y}}{2\sqrt{\pi y}}, \qquad g''(y) = \frac{e^{-y}}{4y\sqrt{\pi y}}(2y + 1) \geq 0.$$

Therefore, we see that

$$(4.17) \qquad \delta g'(\delta) = \mathcal{O}(\delta^{1/2}), \qquad \delta \int_\delta^\infty |g''(y)| \, dy = -\delta g'(\delta) = \mathcal{O}(\delta^{1/2}).$$

Furthermore,

$$
\begin{aligned}
\int_1^\infty g(\delta x) \, dx &= \delta^{-1} \int_\delta^\infty g(y) \, dy = \delta^{-1} \int_\delta^\infty \left( \frac{1}{\sqrt{\pi}} \int_{\sqrt{y}}^\infty e^{-u^2} \, du \right) dy \\
(4.18) \\
&= \delta^{-1} \int_0^\infty \left( \frac{1}{\sqrt{\pi}} \int_{\sqrt{y}}^\infty e^{-u^2} \, du \right) dy - \delta^{-1} \int_0^\delta \left( \frac{1}{\sqrt{\pi}} \int_{\sqrt{y}}^\infty e^{-u^2} \, du \right) dy.
\end{aligned}
$$

Then from $g(\delta) = \tfrac{1}{2} + \mathcal{O}(\delta^{1/2})$, we get

$$(4.19) \qquad \delta^{-1} \int_0^\delta \left( \frac{1}{\sqrt{\pi}} \int_{\sqrt{y}}^\infty e^{-u^2} \, du \right) dy = \frac{1}{2} + \mathcal{O}(\delta^{1/2})$$

and

$$
\begin{aligned}
& \int_0^\infty \left( \frac{1}{\sqrt{\pi}} \int_{\sqrt{y}}^\infty e^{-u^2} \, du \right) dy \\
(4.20) \\
&= \frac{1}{\sqrt{\pi}} y \int_{\sqrt{y}}^\infty e^{-u^2} \, du \Big|_0^\infty - \int_0^\infty y \frac{1}{\sqrt{\pi}} \tfrac{1}{2} y^{-1/2} \cdot -e^{-y} \, dy \\
&= \frac{1}{2\sqrt{\pi}} \int_0^\infty y^{1/2} e^{-y} \, dy = \frac{1}{4}.
\end{aligned}
$$

Therefore,

$$(4.21) \qquad \int_1^\infty g(\delta x) \, dx = \frac{1}{4\delta} - \frac{1}{2} + \mathcal{O}(\delta^{1/2}), \qquad \delta \downarrow 0.$$

It finally follows that

$$(4.22) \quad G(\varepsilon) = \left( \frac{1}{4\delta} - \frac{1}{2} + \mathcal{O}(\delta^{1/2}) \right) + \frac{1}{2} \left( \frac{1}{2} + \mathcal{O}(\delta^{1/2}) \right); \qquad \delta = \frac{1}{2} \varepsilon^2,$$

and we obtain (4.9). It is thus concluded that

$$(4.23) \qquad G(\beta) = \frac{1}{2\beta^2} - \frac{1}{4} - \frac{1}{\sqrt{2\pi}} \sum_{r=0}^\infty \frac{\zeta(-1/2 - r)(-1/2)^r}{r!(2r + 1)} \beta^{2r+1}.$$

Combining (4.1), (4.2) and (4.23), we then obtain

$$
\begin{aligned}
\mathbb{E}M &= \sum_{n=1}^{\infty} \frac{e^{-(1/2)n\beta^2}}{\sqrt{2\pi n}} - \beta \sum_{n=1}^{\infty} P(-\beta\sqrt{n}) \\
&= \sum_{n=1}^{\infty} \frac{e^{-(1/2)n\beta^2}}{\sqrt{2\pi n}} - \beta G(\beta) \\
&= \frac{1}{\beta} + \frac{1}{\sqrt{2\pi}} \sum_{r=0}^{\infty} \frac{\zeta(1/2-r)(-1/2)^r}{r!} \beta^{2r} \\
&\quad - \beta \left[ \frac{1}{2\beta^2} - \frac{1}{4} - \frac{1}{\sqrt{2\pi}} \sum_{r=0}^{\infty} \frac{\zeta(-1/2-r)(-1/2)^r}{r!(2r+1)} \beta^{2r+1} \right] \\
&= \frac{1}{2\beta} + \frac{1}{4}\beta + \frac{1}{\sqrt{2\pi}} \left\{ \sum_{r=0}^{\infty} \frac{\zeta(1/2-r)(-1/2)^r}{r!} \beta^{2r} \right. \\
&\quad \left. + \sum_{r=0}^{\infty} \frac{\zeta(-1/2-r)(-1/2)^r}{r!(2r+1)} \beta^{2r+2} \right\}.
\end{aligned}
\tag{4.24}
$$

Splitting off the term with $r = 0$ and replacing the summation index $r = 1, 2, \ldots$ by $r + 1$, $r = 0, 1, \ldots$ in the first series in (4.24), we get

$$
\mathbb{E}M = \frac{1}{2\beta} + \frac{\zeta(1/2)}{\sqrt{2\pi}} + \frac{1}{4}\beta + \frac{1}{\sqrt{2\pi}} \sum_{r=0}^{\infty} \frac{\zeta(-1/2-r)(-1/2)^r \beta^{2r+2}}{r!(2r+1)(2r+2)}.
\tag{4.25}
$$

## 5. Proof of Theorem 3.

From Spitzer's identity [24], we get

$$
\operatorname{Var} M = \sum_{n=1}^{\infty} \frac{1}{n} \mathbb{E}((S_n^+)^2),
\tag{5.1}
$$

which, using the normality of $S_n$, yields

$$
\begin{aligned}
\operatorname{Var} M &= \sum_{n=1}^{\infty} \frac{1}{n\sqrt{2\pi}} \int_{\beta\sqrt{n}}^{\infty} (x\sqrt{n} - \beta n)^2 e^{-x^2/2} \, dx \\
&= \sum_{n=1}^{\infty} \left( (\beta^2 n + 1) P(-\beta\sqrt{n}) - \frac{\beta}{\sqrt{2\pi}} \sqrt{n} e^{-\beta^2 n/2} \right),
\end{aligned}
\tag{5.2}
$$

where the second equality in (5.2) follows from partial integration. We have established earlier [see (4.6)] that

$$
\frac{1}{\sqrt{2\pi}} \sum_{n=1}^{\infty} \sqrt{n} e^{-\beta^2 n/2} = \frac{1}{\sqrt{2\pi}} \sum_{r=0}^{\infty} \frac{\zeta(-1/2-r)(-1/2)^r}{r!} \beta^{2r} + \frac{1}{\beta^3}.
\tag{5.3}
$$

Therefore, it remains to evaluate

$$(5.4) \qquad I(\beta) = \sum_{n=1}^{\infty} nP(-\beta\sqrt{n}) = \sum_{n=1}^{\infty} \frac{n}{\sqrt{2\pi}} \int_{\beta\sqrt{n}}^{\infty} e^{-x^2/2}\,dx$$

and to combine the results with (5.3) and (4.23) according to (5.2).

There holds

$$(5.5) \qquad \begin{aligned} I'(\beta) &= -\sum_{n=1}^{\infty} \frac{n^{3/2}}{\sqrt{2\pi}} e^{-(1/2)\beta^2 n} \\ &= \frac{-e^{-(1/2)\beta^2}}{\sqrt{2\pi}} \Phi\Big(z = e^{-(1/2)\beta^2}, s = -\frac{3}{2}, v = 1\Big), \end{aligned}$$

and by Bateman's result (2.5),

$$(5.6) \qquad \begin{aligned} I'(\beta) &= \frac{-e^{-(1/2)\beta^2}}{\sqrt{2\pi}} \left[ \frac{\Gamma(5/2)}{e^{-(1/2)\beta^2}} \Big(\frac{1}{2}\beta^2\Big)^{-5/2} \right. \\ &\qquad\qquad \left. + e^{(1/2)\beta^2} \sum_{r=0}^{\infty} \zeta\Big(-\frac{3}{2}-r\Big) \frac{(-(1/2)\beta^2)^r}{r!} \right] \\ &= \frac{-3}{4\sqrt{2}} \Big(\frac{1}{2}\beta^2\Big)^{-5/2} - \frac{1}{\sqrt{2\pi}} \sum_{r=0}^{\infty} \zeta\Big(-\frac{3}{2}-r\Big) \frac{(-(1/2)\beta^2)^r}{r!} \\ &= -3\beta^{-5} - \frac{1}{\sqrt{2\pi}} \sum_{r=0}^{\infty} \frac{\zeta(-3/2-r)(-1/2)^r}{r!} \beta^{2r}, \end{aligned}$$

assuming that $0 < \beta < 2\sqrt{\pi}$. The series on the last line of (5.6) is well behaved in $0 \le \beta < 2\sqrt{\pi}$, whence $I'(\beta) + 3\beta^{-5}$ is integrable, and we obtain

$$(5.7) \qquad \begin{aligned} I(\beta) - \frac{3}{4}\beta^{-4} &= \lim_{\varepsilon\downarrow 0} \Big( I(\varepsilon) - \frac{3}{4}\varepsilon^{-4} \Big) \\ &\qquad - \frac{1}{\sqrt{2\pi}} \sum_{r=0}^{\infty} \frac{\zeta(-3/2-r)(-1/2)^r}{r!(2r+1)} \beta^{2r+1}. \end{aligned}$$

We shall show that

$$(5.8) \qquad \lim_{\varepsilon\downarrow 0}(I(\varepsilon) - \tfrac{3}{4}\varepsilon^{-4}) = -\frac{1}{24}$$

by applying the Euler–Maclaurin summation formula (4.10) with $m = 1$, $N \to \infty$ as before. We consider now

$$(5.9) \qquad f_\delta(x) = \frac{\delta x}{\sqrt{\pi}} \int_{\sqrt{\delta x}}^{\infty} e^{-u^2}\,du =: h(\delta x); \qquad \delta = \frac{1}{2}\varepsilon^2,$$

in which

$$(5.10) \qquad h(x) = xg(x); \qquad g(x) = \frac{1}{\sqrt{\pi}} \int_{\sqrt{x}}^{\infty} e^{-u^2} \, du, \qquad x \geq 0.$$

Then

$$(5.11) \qquad \begin{aligned} I(\varepsilon) = \frac{1}{\delta} \Big[ &\int_1^\infty h(\delta x) \, dx + \frac{1}{2} h(\delta) \\ &- \frac{1}{2} B_2 h'(\delta) \delta - \int_1^\infty \delta^2 h''(\delta x) \frac{B_2(x - \lfloor x \rfloor)}{2} \, dx \Big]. \end{aligned}$$

Next we shall take $\delta \downarrow 0$, and to that end, we see that

$$(5.12) \qquad \frac{1}{\delta} h(\delta) = g(\delta) \to \frac{1}{2}; \qquad h'(\delta) = g(\delta) - \frac{\delta^{1/2}}{2\sqrt{\pi}} e^{-\delta} \to \frac{1}{2}, \qquad \delta \downarrow 0.$$

Furthermore,

$$(5.13) \qquad \frac{1}{\delta} \int_1^\infty h(\delta x) \, dx = \frac{1}{\delta^2} \int_0^\infty h(x) \, dx - \frac{1}{\delta^2} \int_0^\delta h(x) \, dx,$$

in which

$$(5.14) \qquad \frac{1}{\delta^2} \int_0^\delta x g(x) \, dx \to \frac{1}{4} + \mathcal{O}(\delta^{1/2}).$$

Also, by partial integration,

$$(5.15) \qquad \begin{aligned} \int_0^\infty h(x) \, dx &= \int_0^\infty x \frac{1}{\sqrt{\pi}} \Big( \int_{\sqrt{x}}^\infty e^{-u^2} \, du \Big) \, dx \\ &= \frac{x^2}{2\sqrt{\pi}} \int_{\sqrt{x}}^\infty e^{-u^2} \, du \Big|_0^\infty - \int_0^\infty \frac{x^2}{2\sqrt{\pi}} \cdot -\frac{1}{2} x^{-1/2} e^{-x} \, dx \\ &= \frac{1}{4\sqrt{\pi}} \int_0^\infty x^{3/2} e^{-x} \, dx = \frac{1}{4\sqrt{\pi}} \Gamma(5/2) = \frac{3}{16}. \end{aligned}$$

Therefore,

$$(5.16) \qquad \frac{1}{\delta} \int_1^\infty h(\delta x) \, dx = \frac{3}{16\delta^2} - \frac{1}{4} + \mathcal{O}(\delta^{1/2}).$$

Finally,

$$(5.17) \qquad \begin{aligned} h''(x) = (xg(x))'' &= 2g'(x) + xg''(x) \\ &= \frac{1}{2\sqrt{\pi x}} \Big( x - \frac{3}{2} \Big) e^{-x} \in L^1([0, \infty)) \end{aligned}$$

and

$$(5.18) \qquad \frac{1}{2} B_2(x - \lfloor x \rfloor) = \frac{1}{2\pi^2} \sum_{k=1}^\infty \frac{\cos 2\pi kx}{k^2};$$

see [7], page 41. Therefore,

$$
\begin{aligned}
\delta \int_1^\infty & h''(\delta x) \tfrac{1}{2} B_2(x - \lfloor x \rfloor)\, dx \\
(5.19) \qquad & = \int_\delta^\infty h''(x) \tfrac{1}{2} B_2(x/\delta - \lfloor x/\delta \rfloor)\, dx \\
& = \frac{1}{2\pi^2} \sum_{k=1}^\infty \frac{1}{k^2} \int_\delta^\infty h''(x) \cos(2\pi k x/\delta)\, dx \to 0, \qquad \delta \downarrow 0,
\end{aligned}
$$

since $\int_\delta^\infty h''(x) \cos(2\pi k x/\delta)\, dx \to 0$ as $\delta \downarrow 0$ by the Riemann–Lebesgue lemma on Fourier integrals. Putting this altogether, we find (recall $\delta = \frac{1}{2}\varepsilon^2$)

$$
(5.20) \qquad \lim_{\varepsilon \downarrow 0}(I(\varepsilon) - \tfrac{3}{4}\varepsilon^{-4}) = -\tfrac{1}{4} + \tfrac{1}{2} \cdot \tfrac{1}{2} - \tfrac{1}{2} \cdot \tfrac{1}{6} \cdot \tfrac{1}{2} - 0 \ = -\tfrac{1}{24}.
$$

Hence, we obtain, for $0 < \beta < 2\sqrt{\pi}$,

$$
(5.21) \qquad I(\beta) = \frac{3}{4}\beta^{-4} - \frac{1}{24} - \frac{1}{\sqrt{2\pi}} \sum_{r=0}^\infty \frac{\zeta(-3/2 - r)(-1/2)^r}{r!(2r+1)} \beta^{2r+1}.
$$

We insert this result, together with (5.3) and (4.23), into (5.2) and get

$$
\begin{aligned}
\operatorname{Var} M = & -\frac{1}{\sqrt{2\pi}} \sum_{r=0}^\infty \frac{\zeta(-3/2 - r)(-1/2)^r}{r!(2r+1)} \beta^{2r+3} + \frac{3}{4}\beta^{-2} - \frac{\beta^2}{24} \\
(5.22) \qquad & -\frac{1}{\sqrt{2\pi}} \sum_{r=0}^\infty \frac{\zeta(-1/2 - r)(-1/2)^r}{r!(2r+1)} \beta^{2r+1} + \frac{1}{2}\beta^{-2} - \frac{1}{4} \\
& -\frac{1}{\sqrt{2\pi}} \sum_{r=0}^\infty \frac{\zeta(-1/2 - r)(-1/2)^r}{r!} \beta^{2r+1} - \beta^{-2}.
\end{aligned}
$$

Splitting off the terms with $r = 0$ and replacing the summation index $r = 1, 2, \ldots$ by $r + 1$, $r = 0, 1, \ldots$ in the last two series in the right-hand side of (5.22), we get

$$
\begin{aligned}
\operatorname{Var} M = & \frac{1}{4\beta^2} - \frac{1}{4} - \frac{2\zeta(-1/2)}{\sqrt{2\pi}}\beta - \frac{1}{24}\beta^2 \\
& -\frac{1}{\sqrt{2\pi}} \sum_{r=0}^\infty \Bigg\{ \frac{\zeta(-3/2 - r)(-1/2)^r}{r!(2r+1)} \beta^{2r+3} \\
& \qquad\qquad + \frac{\zeta(-3/2 - r)(-1/2)^{r+1}}{(r+1)!(2r+3)} \beta^{2r+3} \\
(5.23) \qquad & \qquad\qquad + \frac{\zeta(-3/2 - r)(-1/2)^{r+1}}{(r+1)!} \beta^{2r+3} \Bigg\}
\end{aligned}
$$

$$= \frac{1}{4\beta^2} - \frac{1}{4} - \frac{2\zeta(-1/2)}{\sqrt{2\pi}}\beta - \frac{1}{24}\beta^2$$

$$- \frac{2}{\sqrt{2\pi}} \sum_{r=0}^{\infty} \frac{\zeta(-3/2-r)(-1/2)^r \beta^{2r+3}}{r!(2r+1)(2r+2)(2r+3)}.$$

**6. Convergence comparison Spitzer formulas and Lerch series.** It is immediately clear that the infinite series in (2.1), (2.2) and (2.3) converge more rapidly for smaller values of $\beta$, while the contrary holds for their Spitzer-type counterparts (3.1), (4.1) and (5.1). To exemplify this difference in speed of convergence, we consider (4.2), that is,

$$(6.1) \qquad \sum_{n=1}^{\infty} \frac{e^{-(1/2)n\beta^2}}{\sqrt{2\pi n}} = \frac{1}{\sqrt{2\pi}} \sum_{r=0}^{\infty} \frac{\zeta(1/2-r)}{r!}\left(-\frac{\beta^2}{2}\right)^r.$$

The left-hand side series converges for all $\beta > 0$, while the right-hand side series converges for all $\beta \in \mathbb{C}$, $|\beta| < 2\sqrt{\pi}$. From [27] Section 13.151 (page 269),

$$(6.2) \qquad 2^{1-s}\Gamma(s)\zeta(s)\cos(\tfrac{1}{2}s\pi) = \pi^s\zeta(1-s).$$

With $s = r + \frac{1}{2}$, the asymptotics of the $\Gamma$-function and the fact that $\zeta(r + 1/2) \to 1$ as $r \to \infty$, we see that

$$(6.3) \qquad \left|\frac{1}{\sqrt{2\pi}}\frac{\zeta(1/2-r)}{r!}\left(-\frac{\beta^2}{2}\right)^r\right| \approx \frac{1}{\pi\sqrt{2r+1}}\left(\frac{\beta^2}{4\pi}\right)^r, \qquad r \to \infty.$$

Hence, for comparing the convergence rates of the two series in (6.1), it is enough to find the point $\beta_0 > 0$ such that

$$(6.4) \qquad e^{-(1/2)\beta_0^2} = \frac{\beta_0^2}{4\pi}.$$

With $x = \frac{1}{2}\beta^2$, we need to solve $x_0 e^{x_0} = 2\pi$ with $x_0 > 0$. This yields $x_0 = 1.4597$, $\beta_0 = 1.7086$, and the common value of the two members in (6.4) equals $0.2323$. See [5], Section 2, where a similar strategy is developed in connection with the evaluation of Legendre's chi-function.

**7. Conclusions and outlook.** We have presented analytic formulas of the Chang–Peres type, involving the Riemann zeta function, for the quantities

$$(7.1) \qquad J_k(\beta) = \sum_{n=1}^{\infty} \frac{1}{n}\mathbb{E}((S_n^+)^k), \qquad k = 0, 1, 2,$$

yielding $\mathbb{P}(M = 0)$, $\mathbb{E}M$ and $\operatorname{Var} M$ of the maximum $M$ of the standard Gaussian random walk with negative drift $-\beta$. The quantities can be expressed using the normality of $S_n$ as

$$(7.2) \qquad J_k(\beta) = \sum_{n=1}^{\infty} \frac{n^{k-1/2}}{\sqrt{2\pi}} \int_{\beta}^{\infty} (y-\beta)^k e^{-(1/2)ny^2}\, dy.$$

The following general result can be shown: There holds for $k = 1, 2, \dots$

$$
\begin{aligned}
J_k(\beta) = {} & \frac{(k-1)!}{(2\beta)^k} \\
(7.3) \quad & + \sum_{j=0}^{k} \binom{k}{j} \frac{(-1)^j \Gamma((k-j+1)/2)}{\sqrt{2\pi}} \zeta\left(-\frac{1}{2}k - \frac{1}{2}j + 1\right) 2^{(k-j-1)/2} \beta^j \\
& + \frac{(-1)^{k+1} k!}{\sqrt{2\pi}} \sum_{r=0}^{\infty} \frac{\zeta(-k-r+1/2)(-1/2)^r \beta^{2r+k+1}}{r!(2r+1)\cdots(2r+k+1)},
\end{aligned}
$$

when $0 < \beta < 2\sqrt{\pi}$. That such a formula should hold for $J_k(\beta)$ follows from differentiation of (7.2) $k+1$ times so that there results

$$
(7.4) \qquad J_k^{(k+1)}(\beta) = (-1)^{k+1} k! \sum_{n=1}^{\infty} \frac{n^{k-1/2}}{\sqrt{2\pi}} e^{-(1/2)n\beta^2},
$$

in which the right-hand side is readily expressed in terms of Lerch's transcendent as we did before. Then, using (2.5) and integrating the identity thus obtained $k+1$ times, we arrive at (7.3), except for the $k+1$ integration constants that appear at the right-hand side of (7.3) as the coefficients of $\beta^j$, $j = 0, 1, \dots, k$. The actual determination of these integration constants is still based on the Euler–Maclaurin summation formula in its general form (4.10), but is so complicated that a full proof of (7.3) is outside the scope of the present paper.

We are presently undertaking an effort to analyze a specific queueing model under a heavy traffic scaling. In the language of Queueing Theory we are dealing with a discrete-time bulk service queue with batch arrivals in which the arrival process is a Poisson process whose arrival rate $\lambda$ is just slightly smaller than the service capacity $s$. In the case that $s = \lambda + \beta\sqrt{\lambda}$, with $\beta > 0$ fixed and $\lambda \to \infty$ (Halfin–Whitt regime, see [16, 18]), the equilibrium distribution of the queue converges to that of the Gaussian random walk. The analysis of this equilibrium distribution for finite $\lambda$ is, however, far more complicated than in the case of the Gaussian random walk. Already the first two moments require an evaluation in terms of Lerch's transcendent of the expressions

$$
(7.5) \qquad\qquad \sum_{n=1}^{\infty} n^{k-1/2} \int_{\beta}^{\infty} e^{-(1/2)nx^2} x^i \, dx,
$$

for all integer $k$ and all $i = 0, 1, \dots$. The results obtained by us for the expressions (7.5) allow us to fully establish (7.3), that is, including the integration constants.

A further comment concerns the restricted validity of formulas (2.1)–(2.3), namely for $0 < \beta < 2\sqrt{\pi}$, while the corresponding Spitzer formulas (3.1), (4.1) and (5.1) make sense for all $\beta > 0$. By using (6.2), the explicit formula $\zeta(s) = \sum_{n=1}^{\infty} n^{-s}$ valid for $s > 1$, and some further manipulations, the infinite series occurring in (2.1)–(2.3) can be re-expressed as follows. We have for $\beta > 0$

$$
\begin{aligned}
(7.6) \quad & \frac{\beta}{\sqrt{2\pi}} \sum_{r=0}^{\infty} \frac{\zeta(1/2-r)(-(1/2)\beta^2)^r}{r!(2r+1)} \\
& = \frac{\zeta(1/2)}{\sqrt{2\pi}}\beta + \frac{\beta}{\pi} \operatorname{Re}\left[e^{\pi i/4} S_0\left(\frac{-i\beta^2}{4\pi}\right)\right],
\end{aligned}
$$

$$
\begin{aligned}
(7.7) \quad & \frac{\beta^2}{\sqrt{2\pi}} \sum_{r=0}^{\infty} \frac{\zeta(-1/2-r)(-(1/2)\beta^2)^r}{r!(2r+1)(2r+2)} \\
& = \frac{\beta^2}{2\pi^2} \operatorname{Re}\left[-e^{\pi i/4} S_1\left(\frac{-i\beta^2}{4\pi}\right)\right]
\end{aligned}
$$

and

$$
\begin{aligned}
(7.8) \quad & \frac{2\beta^3}{\sqrt{2\pi}} \sum_{r=0}^{\infty} \frac{\zeta(-3/2-r)(-(1/2)\beta^2)^r}{r!(2r+1)(2r+2)(2r+3)} \\
& = \frac{\beta^3}{4\pi^3} \operatorname{Re}\left[e^{\pi i/4} S_2\left(\frac{-i\beta^2}{4\pi}\right)\right],
\end{aligned}
$$

in which

$$
(7.9) \quad S_0(b) = \frac{\sqrt{\pi}}{\sqrt{b}} \sum_{n=1}^{\infty} (\arcsin(b/n)^{1/2} - (b/n)^{1/2}),
$$

$$
(7.10) \quad S_1(b) = \frac{\sqrt{\pi}}{2b} \sum_{n=1}^{\infty} \frac{1}{n}(\sqrt{n} - \sqrt{n-b}),
$$

$$
(7.11) \quad S_2(b) = \frac{\sqrt{\pi}}{4b} \sum_{n=1}^{\infty} \frac{1}{n^2}(\sqrt{n} - \sqrt{n-b})
$$

are well defined functions for $b$ ($= -i\beta^2/4\pi$) in an open set containing the imaginary axis. (These alternative expressions are intimately related with Lerch's transformation formula ([13], 1.11(7) on page 29) and complement the results obtained through analytic continuation by Chang and Peres [10].) Thus, for values of $\beta$ larger than $2\sqrt{\pi}$, we can evaluate $\mathbb{P}(M = 0)$, $\mathbb{E}M$, $\operatorname{Var}M$ using (2.1)–(2.3) and (7.6)–(7.11). Although the series in (7.9)–(7.11) converge slowly, they can be evaluated quite conveniently by using dedicated forms of the Euler–Maclaurin summation formula (4.10). We may also note that the infinite series occurring in the expression (7.3) for the general $J_k(\beta)$ can be re-expressed in a similar fashion.

Finally, let us return once more to the rather special constant $-\zeta(1/2)/\sqrt{2\pi} \approx 0.5826$ (see Section 1). From Theorem 2 [or even (1.2)], we get

$$(7.12) \qquad \frac{1}{2\beta} - \mathbb{E}M = -\frac{\zeta(1/2)}{\sqrt{2\pi}} + \mathcal{O}(\beta),$$

where $1/(2\beta)$ is known to be an upper bound on $\mathbb{E}M$ that becomes tight in heavy traffic, that is, for $\beta \downarrow 0$ (see [19, 20]). Hence, the difference between the upper bound and the true value tends to a nontrivial constant $-\zeta(1/2)/\sqrt{2\pi}$. In [12], constants of this type are studied for a general class of random walks (possibly non-Gaussian) with zero drift. From Theorem 3 and the general result (7.3), we see that such nontrivial constants exist for higher moments of $M$ as well.

**Acknowledgment.** We are grateful to Bert Zwart for discussions that motivated this paper and for pointing us to the references [10, 18].

Digital Signal Processing Group
Philips Research WO-02
5656 AA Eindhoven
The Netherlands
E-mail: a.j.e.m.janssen@philips.com

Eurandom
5612 AZ Eindhoven
The Netherlands
E-mail: j.s.h.v.leeuwaarden@tue.nl